\def\titlep{Unitary isomorphism of Fock spaces of bosons and fermions
arising from a representation of the Cuntz algebra $\co{2}$}
\documentclass[11pt]{article}
\usepackage{graphicx,ifthen}
\usepackage{amssymb}

\font\germ=eufm10 at12pt

\def\goth#1{\hbox{\germ#1}}

\newcommand{\sdag}{\scriptsize \dag}
\newcommand{\qed}{\hbox{\rule[-2pt]{3pt}{6pt}}}
\newcommand{\qedh}{\hfill\qed \\}

\newtheorem{Thm}{Theorem}[section]

\newtheorem{rem}[Thm]{Remark}

\newtheorem{ex}[Thm]{Example}
\newtheorem{defi}[Thm]{Definition}
\newtheorem{lem}[Thm]{Lemma}

\newtheorem{prop}[Thm]{Proposition}
\newtheorem{prob}[Thm]{Problem}
\newtheorem{cor}[Thm]{Corollary}

\def\cal#1{\mathcal #1}
\def\con{{\cal O}_{N}}
\def\coni{{\cal O}_{\infty}}
\def\edot{=1,\ldots,N}
\def\pr{{\it Proof.}\quad}

\def\scm#1{S({\bf C}^{N})^{\otimes #1}}

\def\co#1{{\cal O}_{#1}}
\def\ltn{l_{2}({\bf N})}

\def\disp#1{{\displaystyle #1}}
\def\brl{branching law}
\def\bfsnl{{\rm BFS}_{N}(\Lambda)}

\addtocounter{footnote}{1}
\def\sftt#1{
\setcounter{equation}{0}
\addtocounter{footnote}{1}
\section{#1}
}

\def\ssft#1{\subsection{#1}}
\def\sssft#1{\subsubsection{#1}}

\begin{document}
%
% Personal data
%
\def\autherp{Katsunori Kawamura}
\def\emailp{Electronic mail: kawamura@kurims.kyoto-u.ac.jp.}
\def\addressp{{\it {\small College of Science and Engineering Ritsumeikan University,}}\\
{\it {\small 1-1-1 Noji Higashi, Kusatsu, Shiga 525-8577, Japan}}
}

\def\infw{\Lambda^{\frac{\infty}{2}}V}
\def\zhalfs{{\bf Z}+\frac{1}{2}}
\def\ems{\emptyset}
\def\pmvac{|{\rm vac}\!\!>\!\! _{\pm}}
\def\vac{|{\rm vac}\rangle _{+}}
\def\dvac{|{\rm vac}\rangle _{-}}
\def\ovac{|0\rangle}
\def\tovac{|\tilde{0}\rangle}
\def\expt#1{\langle #1\rangle}
\def\zph{{\bf Z}_{+/2}}
\def\zmh{{\bf Z}_{-/2}}
\def\brl{branching law}
\def\bfsnl{{\rm BFS}_{N}(\Lambda)}
\def\scm#1{S({\bf C}^{N})^{\otimes #1}}
\def\mqb{\{(M_{i},q_{i},B_{i})\}_{i=1}^{N}}
\def\zhalf{\mbox{${\bf Z}+\frac{1}{2}$}}
\def\zmha{\mbox{${\bf Z}_{\leq 0}-\frac{1}{2}$}}
\newcommand{\mline}{\noindent
\thicklines
\setlength{\unitlength}{.1mm}
\begin{picture}(1000,5)
\put(0,0){\line(1,0){1250}}
\end{picture}
\par
 }
\def\sd#1{#1^{\sdag}}
\def\dlim{{\cal D}\mbox{-}\lim}
\def\dsum{{\cal D}\mbox{-}\sum}
\def\fsum{F\mbox{-}\sum}
\def\bx{\mbox{\boldmath$x$}}

%
%%%%%%%%% Cut from here %%%%%%%%%%
%\input comm.txt
%%%%%%%%% End of Cut %%%%%%%%%
%
%
\setcounter{section}{0}
\setcounter{footnote}{0}
\setcounter{page}{1}
\pagestyle{plain}

%%%%%%%%%%%%%%%%%%%%%%%%%%%%%%%%%%%%%%%%%%%%%%%%%%%%%%%%%%%
%
% Title
%
\title{\titlep}
\author{\autherp\thanks{\emailp}
\\ 
\addressp
}
\date{}
\maketitle

%%%%%%%%%%%%%%%%%%%%%%%%%%%%%%%%%%%%%%%%%%%%%%%%%%%%%%%%%%%
%
% Abstract
%
\begin{abstract}
Bosons and fermions are described 
by using canonical generators of Cuntz algebras
on any permutative representation.
According to branching laws associated with these descriptions,
a certain representation of the Cuntz algebra $\co{2}$ induces
Fock representations ${\cal H}_{B}$ and ${\cal H}_{F}$
of bosons and fermions simultaneously.
From this,
a unitary operator $U$ from ${\cal H}_{B}$ to ${\cal H}_{F}$
is obtained.
We show the explicit formula of the action of $U$
on the standard basis of ${\cal H}_{B}$. 
It is shown that $U$ preserves the particle number
of ${\cal H}_{B}$ and ${\cal H}_{F}$.
\end{abstract}


%%%%%%%%%%%%%%%%%%%%%%%%%%%%%%%%%%%%%%%%%%%%%%%%%%%%%%%%%%%%%
%
% Section 1
%
\sftt{Introduction}
\label{section:first}
We have studied relations among representations
of Cuntz algebras, bosons and fermions 
from a standpoint of branching law.
Even to this day,
algebras of bosons and fermions
are most fundamental ingredients in quantum field theory \cite{NO,PS}.
Their Fock representations are basic and essential in several models.
Therefore it is expected that
a mathematical study of these representations give
a new point of view to quantum field theory.
In this paper,
we show their relations by using the representation theory 
of Cuntz algebras.
A unitary operator from the Bose-Fock space to the Fermi-Fock space
is constructed by using a certain representation of the Cuntz algebra $\co{2}$
which preserves particle number.

%%%%%%%%%%%%%%%%%%%%%%%%%%%%%%%%%%%%%%%%%%%%%%%%%%%%%%%%%%%%%%%
%
% subsection 1.1
%
\ssft{Motivation}
\label{subsection:firstone}
We explain our motivation without mathematical definitions
in this subsection.
Rigorous definitions
will be given in $\S$ \ref{subsection:firsttwo} 
and $\S$ \ref{section:second}.

%%%%%%%%%%%%%%%%%%%%%%%%%%%%%%%%%%%%%%%%%%%%%%%%%%%%%%%%%%%%%
%
% subsection 1.1.1
%
\sssft{Fermions and bosons described by Cuntz algebras}
\label{subsubsection:firstoneone}
In \cite{AK1} and \cite{RBS01}, we have described fermions and bosons 
by using Cuntz algebras
$\co{2}$ and $\coni$, respectively.
Fermions are described as polynomials in canonical generators of $\co{2}$
and their conjugates.
Hence the algebra of fermions is embedded 
into $\co{2}$ as a unital $*$-subalgebra.
Bosons are described as formal power series in
canonical generators of $\coni$ and their conjugates.
This does not mean that the algebra of bosons is embedded into $\coni$.
However, these formal sums make sense 
on a certain dense subspace of any permutative representations of $\coni$.
These two descriptions enable studies of restrictions of representations of
$\co{2}$ and $\coni$ on fermions and bosons.
We have shown that their
Fock representations are derived from certain representations
of these Cuntz algebras.

On the other hand, $\coni$ is embedded into $\co{2}$ 
by using a certain unital $*$-embedding
such that canonical generators of $\coni$ are written as monomials 
in those of $\co{2}$.
With respect to this embedding,
the restriction of any permutative representation of $\co{2}$
is also a permutative representation of $\coni$.
This implies that any permutative representation of $\co{2}$
can be restricted on bosons also.
From these, both fermions and bosons are 
described by using canonical generators of $\co{2}$.
Hence we have an extreme interest in
the comparison between fermions and bosons which 
are simultaneously represented 
on a given representation space of $\co{2}$.

%%%%%%%%%%%%%%%%%%%%%%%%%%%%%%%%%%%%%%%%%%%%%%%%%%%%%%%%%%%%%
%
% subsection 1.1.2
%
\sssft{Bosonization, fermionization and boson-fermion correspondence}
\label{subsubsection:firstonetwo}
In physics,
a description of a given (algebraic) system 
by using bosons ({\it resp.} fermions) is 
called {\it bosonization} \cite{Kopietz,Stone}
({\it resp}. {\it fermionization} \cite{AKMS}).
Especially,
fermions and bosons are often rewritten each other in some senses.
In this case,
a pair of bosonization and fermionization 
is called a {\it boson-fermion correspondence} \cite{Kac,MJD,Oko01}
in the broad sense of the term.

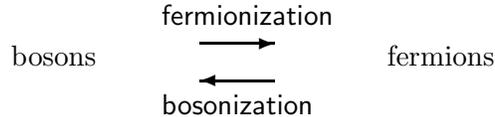
\def\pict{
\put(0,0){bosons}
\put(500,0){fermions}
\put(250,30){\vector(1,0){100}}
\put(200,55){\textsf{fermionization}}
\put(350,-20){\vector(-1,0){100}}
\put(200,-65){\textsf{bosonization}}
}
\thicklines
\setlength{\unitlength}{.1mm}
\begin{picture}(1000,200)(100,10)
\put(300,100){\pict}
\end{picture}

\noindent
These metamorphoses are not only interesting phenomena 
but also important techniques in physics \cite{Coleman,Mandelstam}.
On the other hand,
such descriptions are not well-understood as mathematics.
In almost all cases,
a bosonization (or a fermionization) is done as a computation technique.
Therefore nobody has given its mathematical meaning.
The algebra of bosons and that of fermions are neither isomorphic
nor embedded each other.
Hence neither bosonization nor  fermionization
are executed in a purely algebraic sense without any representation.
Remark that infinite sums and normal orders 
of elements in these algebras
make no sense as elements in algebras. 
These descriptions are usually given on a certain representation space 
with a certain operator topology or as a formal operation \cite{IT,MJD,Oko01}.
Hence a mathematical generalization of these is stated as follows:
%
% Problem 1.1
%
\begin{prob}
\label{prob:first}
Let $A_{1}$ and $A_{2}$ be algebras such that
there is no embedding of one to the other.
Let $V$ be a vector space
and let $\pi_{i}$ be a representation of $A_{i}$ 
on $V$ for $i=1,2$.
\begin{enumerate}
%(i)
\item
For $x\in A_{1}$,
write $\pi_{1}(x)$ by using elements in $\pi_{2}(A_{2})$
(and some exclusive operation)
as a rigorous mathematical statement.  
%(ii)
\item
If (i) is done, then
characterize such description in some senses.
\end{enumerate}
\end{prob}

\noindent
If $\pi_{2}$ is irreducible,
then we can always give an answer to
Problem \ref{prob:first}(i) by using the strong operator topological limit
(\cite{RBS01}, $\S$ I.A.2) when $V$ is a Hilbert space.

In order to avoid difficulties 
about operator-valued distributions and physical assumptions, 
we consider a boson-fermion correspondence between the Bose-Fock space
and the Fermi-Fock space as the simplest example in this paper.

%%%%%%%%%%%%%%%%%%%%%%%%%%%%%%%%%%%%%%%%%%%%%%%%%%
%
%  subsection 1.2
%
\ssft{Recursive boson system and recursive fermion system}
\label{subsection:firsttwo}
We briefly explain recursive boson system and recursive fermion system 
in this subsection.
For $2\leq N\leq \infty$, let $\con$ denote the Cuntz algebra \cite{C}.
Let $\{s_{n}:n\in {\bf N}\}$ and 
$\{t_{1},t_{2}\}$ denote the canonical generators
of $\coni$ and $\co{2}$, respectively, that is,
they satisfy that
%
% Equation 1.1 and 1.2
% 
\begin{eqnarray}
\label{eqn:oinfty}
s_{i}^{*}s_{j}=\delta_{ij}I\quad (i,j\in {\bf N}),\quad 
&\disp{\sum_{m=1}^{k}s_{m}s_{m}^{*}\leq I}
&\quad (\mbox{for any }k\geq 1),\\
\nonumber
\\
\label{eqn:otwo}
t_{i}^{*}t_{j}=\delta_{ij}I\quad  (i,j=1,2),\quad & t_{1}t_{1}^{*}+
t_{2}t_{2}^{*}=I,& 
\end{eqnarray}
where ${\bf N}\equiv \{1,2,3,\ldots\}$.
Let $\{b_{n}:n\in {\bf N}\}$ and $\{a_{n}:n\in {\bf N}\}$
denote bosons and fermions, that is,
%
% Equation 1.3 and 1.4
% 
\begin{eqnarray}
\label{eqn:boson}
\,b_{n}b_{m}^{*}-b_{m}^{*}b_{n}=\delta_{nm}I,\quad  &b_{n}b_{m}-b_{m}b_{n}
=b_{n}^{*}b_{m}^{*}-b_{m}^{*}b_{n}^{*}=0,\\
\nonumber
\\
\label{eqn:fermion}
a_{n}a_{m}^{*}+a_{m}^{*}a_{n}=\delta_{nm}I,\quad &
a_{n}a_{m}+a_{m}a_{n}=a_{n}^{*}a_{m}^{*}+a_{m}^{*}a_{n}^{*}=0
\end{eqnarray}
for each $n,m\in {\bf N}$. 
We described $\{b_{n}:n\in {\bf N}\}$ and $\{a_{n}:n\in {\bf N}\}$
by using $\{s_{n}:n\in {\bf N}\}$ and $\{t_{1},t_{2}\}$ 
in \cite{RBS01,AK1}, respectively
as follows:
%
% Equation 1.5, 1.6
%
\begin{eqnarray}
\label{eqn:rbszero}
b_{1}= \sum_{m=1}^{\infty}\sqrt{m}\,s_{m}s_{m+1}^{*},& 
\,\,b_{n}= \rho(b_{n-1})\quad(n\geq 2),\\ \nonumber
\\
\label{eqn:rfszero}
a_{1}= t_{1}t_{2}^{*},\,\,\quad\qquad\qquad & a_{n}=\zeta(a_{n-1})\quad(n\geq 2)
\end{eqnarray}
where 
%
% Equation 1.7 and 1.8
%
\begin{eqnarray}
\label{eqn:rho}
\rho(x)= &\disp{\sum_{m=1}^{\infty}s_{m}xs_{m}^{*}\qquad
(x\in \coni)},\\
\nonumber
\\
\zeta(y)= &t_{1}yt_{1}^{*}-t_{2}yt_{2}^{*}\quad(y\in \co{2}).
\label{eqn:zeta}
\end{eqnarray}
We call these descriptions of $\{b_{n}:n\in {\bf N}\}$ and $\{a_{n}:n\in {\bf N}\}$
as the {\it recursive boson system (=RBS)} and
the {\it recursive fermion system (=RFS)}, respectively. 
%
% Remark 1.2
%
\begin{rem}
\label{rem:two}
{\rm
\begin{enumerate}
%(i)
\item
Remark that $\{b_{n}:n\in {\bf N}\}$ in (\ref{eqn:rbszero}) 
and $\rho(x)$ in (\ref{eqn:rho})
are not well-defined in $\coni$, 
but they make sense as operators on a certain dense subspace 
of any permutative representation of $\coni$ by Fact 1.1 in \cite{RBS01}.
%(ii)
\item
The $*$-algebra ${\cal B}$ of bosons can {\it never} be
embedded into the $*$-algebra ${\cal A}$ of fermions (\cite{RBS01}, $\S$ 1.1.1).
Especially, ${\cal B}$ and ${\cal A}$ are not $*$-isomorphic.
From this,
bosonization and fermionization are usually executed 
as infinite operations on suitable representation
spaces with respect to certain operator topologies.
\end{enumerate}
}
\end{rem}
%%%%%%%%%%%%%%%%%%%%%%%%%%%%%%%%%%%%%%%%%%%%%%%%%%%%%
%
% subsection 1.3
%
\ssft{Fock representations arising from a representation of $\co{2}$}
\label{subsection:firstthree}
In this subsection,
we show relations among the Bose-Fock representation, 
the Fermi-Fock representation
and a certain representation of $\co{2}$.
Let ${\cal B}$ and ${\cal A}$
be as in Remark \ref{rem:two}(ii) and let
$\{t_{1},t_{2}\}$ be as in (\ref{eqn:otwo}).
Assume that $({\cal H},\pi)$ is a $*$-representation
of $\co{2}$ with a cyclic vector $\Omega$ satisfying 
%
% Equation 1.9
%
\begin{equation}
\label{eqn:one}
\pi(t_{1})\Omega=\Omega.
\end{equation}
This representation exists uniquely up to unitary equivalence 
\cite{BJ,DaPi2,DaPi3}.
We will show an example of this in $\S$ \ref{subsection:fourththree}.
In \cite{AK1,RBS01},
we proved the following for $({\cal H},\pi)$ and $\Omega$ in (\ref{eqn:one}):
\begin{enumerate}
%(i)
\item
Define the dense subspace ${\cal D}$ of ${\cal H}$ by
%
% Equation 1.10
%
\begin{equation}
\label{eqn:cald}
{\cal D}\equiv {\rm Lin}\langle
\{t_{i_{1}}\cdots t_{i_{m}}\Omega:i_{1},\ldots,i_{m}=1,2,\,m\geq 1\}
\rangle.
\end{equation}
Let $\{s_{n}:n\in {\bf N}\}$ be as in (\ref{eqn:oinfty}).
Assume that $\coni$ is embedded into $\co{2}$ by
%
% Equation 1.11
%
\begin{equation}
\label{eqn:embedding}
s_{m}= t_{2}^{m-1}t_{1}\quad(m\geq 1)
\end{equation}
where we define $t_{2}^{0}=I$.
From this, we obtain the
restriction $({\cal H},\pi|_{\coni})$ of 
$({\cal H},\pi)$ on $\coni$.
Furthermore,
(\ref{eqn:rbszero}) and (\ref{eqn:rho}) are rewritten as follows:
%
% Equation 1.12, 1.13
%
\begin{eqnarray}
\label{eqn:rbszerotwo}
b_{1}= \sum_{m=1}^{\infty}\sqrt{m}\,t_{2}^{m-1}t_{1}t_{1}^{*}(t_{2}^{*})^{m},
\quad \,\,b_{n}= \rho(b_{n-1})\quad(n\geq 2),\\ 
\nonumber
\\
\label{eqn:rhotwo}
\rho(x)= \disp{\sum_{m=1}^{\infty}t^{m-1}_{2}t_{1}xt_{1}^{*}(t_{2}^{*})^{m-1}}.
\qquad\qquad 
\end{eqnarray}
By using (\ref{eqn:embedding}),
a unital $*$-representation $\pi_{RBS}$ of ${\cal B}$ on ${\cal D}$
is defined by using (\ref{eqn:rbszero}) 
with respect to $\{\pi|_{\coni}(s_{n}):n\in {\bf N}\}$
and their conjugates on ${\cal D}$.
%(ii)
\item
The unital $*$-representation $\pi_{RFS}$ of ${\cal A}$ on ${\cal H}$ 
is defined by using (\ref{eqn:rfszero})
with respect to $\{\pi(t_{1}),\pi(t_{2})\}$ and their conjugates.
%(iii)
\item
Let $({\cal H}_{B},\pi_{B})$ and $({\cal H}_{F},\pi_{F})$ denote 
the Fock representations of bosons and fermions with vacua 
$\Omega_{B}$ and $\Omega_{F}$,
respectively \cite{BR}, that is,
which are cyclic vectors of dense subspaces of 
${\cal H}_{B}$ and ${\cal H}_{F}$, respectively and 
%
% Equation 1.14
%
\begin{equation}
\label{eqn:two}
\pi_{B}(b_{n})\Omega_{B}=0,\quad 
\pi_{F}(a_{n})\Omega_{F}=0\quad\mbox{for all }n\in {\bf N}.
\end{equation}
Then there exist two unitaries
$V_{B}:{\cal H}\to {\cal H}_{B}$ 
and $V_{F}:{\cal H}\to {\cal H}_{F}$ such that
%
% Equation 1.15
%
\begin{equation}
\label{eqn:vbtwo}
V_{B}\pi_{RBS}(\cdot)V_{B}^{*}=\pi_{B},\quad 
V_{F}\pi_{RFS}(\cdot)V_{F}^{*}=\pi_{F},
\end{equation}
that is, two unitary equivalences 
$\pi_{RBS}\cong \pi_{B}$ and $\pi_{RFS}\cong \pi_{F}$ hold,
and $V_{B}\Omega=\Omega_{B}$ and $V_{F}\Omega=\Omega_{F}$
($\S$ 3.3 of \cite{AK1} and Proposition 3.2 of \cite{RBS01}). 
\end{enumerate}
From these results,
we define the unitary $U$ from ${\cal H}_{B}$ to ${\cal H}_{F}$ by

%
% Equation 1.16
%
\begin{equation}
\label{eqn:unitary}
\begin{picture}(1000,130)(0,260)
\put(800,260){$U\equiv V_{F}V_{B}^{*}$}
\put(360,400){${\cal H}$}
\put(70,250){${\cal H}_{B}$}
\put(600,250){${\cal H}_{F}$,}
\put(270,370){\vector(-2,-1){120}}
\put(480,370){\vector(2,-1){120}}
\put(230,260){\vector(1,0){300}}
%\put(230,270){\vector(-1,0){0}}
\put(360,340){$\circlearrowleft$}
\put(120,360){$V_{B}$}
\put(570,360){$V_{F}$}
\put(370,290){U}
\end{picture}
\end{equation}

\noindent
It seems that three triplets 
$({\cal H},\pi,\Omega)$, $({\cal H}_{B},\pi_{B},\Omega_{B})$ and 
$({\cal H}_{F},\pi_{F},\Omega_{F})$ are similar
in a sense of the representation theory of operator algebras
because they are irreducible representations
with cyclic vectors $\Omega,\Omega_{B},\Omega_{F}$ such that they are uniquely determined 
by algebraic equations with respect to $\Omega,\Omega_{B},\Omega_{F}$.
Remark that both $V_{B}$ and $V_{F}$ in (\ref{eqn:vbtwo})
are naturally constructed from (\ref{eqn:two}).
From (\ref{eqn:rbszero}) and (\ref{eqn:rfszero}),
both equations in (\ref{eqn:two}) are derived from (\ref{eqn:one}).
In this sense, (\ref{eqn:one}) is most fundamental
from the perspective of the representation theory of these algebras.

In spite that two unitaries $V_{B}$ and $V_{F}$ have an established role
of unitary equivalence of representations,
we can not explain a meaning of the unitary $U$ in (\ref{eqn:unitary}) clearly.
By computing $U$ on the standard basis of ${\cal H}_{B}$,
we consider {\it how similar 
$({\cal H}_{B},\pi_{B},\Omega_{B})$ and $({\cal H}_{F},\pi_{F},\Omega_{F})$ are}
in the next subsection.

%%%%%%%%%%%%%%%%%%%%%%%%%%%%%%%%%%%%%%%%%%%%%%%%%%%%%%%%%%
%
% subsection 1.4
%
\ssft{Main theorem}
\label{subsection:firstfour}
In this subsection, we show our main theorem.
Let $\pi_{B}$ and $\pi_{F}$ be as in (\ref{eqn:two}).
We identify $\pi_{B}(b_{n})$ and $\pi_{F}(a_{n})$
with $b_{n}$ and $a_{n}$ for $n\geq 1$, respectively. 
%
% Theorem 1.3
%
\begin{Thm}
\label{Thm:main}
\begin{enumerate}
%(i)
\item
The unitary $U$ in (\ref{eqn:unitary}) satisfies that
$U\Omega_{B}=\Omega_{F}$ and 
%
% Equation 1.17
%
\begin{equation}
\label{eqn:block}
\hspace{-.9cm}
U\left(\prod_{i=1}^{m}(b_{n_{i}}^{*})^{k_{i}}\right)\Omega_{B}
=C\cdot A_{n_{1}-1,k_{1}}
A_{n_{2}+k_{1}-1,k_{2}}
\cdots 
A_{n_{m}+k_{1}+\cdots+k_{m-1}-1,k_{m}}\Omega_{F}
\end{equation}
for $1\leq n_{1}<\cdots<n_{m}$ and $k_{1},\ldots,k_{m}\geq 1$
where $C$ denotes the normalization constant which is given by
$C\equiv \sqrt{k_{1}!\cdots k_{m}!}$ and
%
% Equation 1.18
%
\begin{equation}
\label{eqn:blocktwo}
A_{n,m}\equiv a_{n+1}^{*}\cdots a_{n+m}^{*}\quad(n,m\geq 1).
\end{equation}
%
%(ii)
\item
Let ${\cal H}_{B,n}$ and ${\cal H}_{F,n}$
denote the subspaces of all $n$-particle states 
of ${\cal H}_{B}$ and ${\cal H}_{F}$ for $n\geq 0$, respectively.
Then 
%
% Equation 1.19
%
\begin{equation}
\label{eqn:particle}
U{\cal H}_{B,n}={\cal H}_{F,n}\quad(n\geq 0),
\end{equation}
that is, $U$ satisfies the particle number conservation law between
bosons and fermions.
\end{enumerate}
\end{Thm}

\noindent
For other isomorphism theorems of Fock spaces,
see $\S$ 14.10 of \cite{Kac} and Theorem 5.1 in \cite{MJD}. 
From Theorem \ref{Thm:main}, the following holds.
%
% Corollary 1.4
%
\begin{cor}
\label{cor:one}
Up to vacuum vectors $\Omega_{B}$, $\Omega_{F}$ 
and the normalization constant,
monomials of bosons and fermions are one-to-one corresponded 
in the Fock space as follows:
%
% Equation 1.20
%
\begin{equation}
\label{eqn:threetwentyeight}
(b_{n_{1}}^{*})^{k_{1}}\cdots (b_{n_{m}}^{*})^{k_{m}}
\mapsto A_{n_{1}-1,k_{1}}\cdots 
A_{n_{m}+k_{1}+\cdots+k_{m-1}-1,k_{m}}.
\end{equation}
\end{cor}

%
% Remark 1.5
% 
\begin{rem}
\label{rem:one}
{\rm
\begin{enumerate}
%(i)
\item
From Theorem \ref{Thm:main}(i),
we see that a mode $b_{n}^{*}$ of a boson is nearly corresponded with
a block $A_{n^{'},m^{'}}$ of fermions in (\ref{eqn:blocktwo}).
%(ii)
\item
The statement of Theorem \ref{Thm:main}(i) is
a result of the definition of $U$ in (\ref{eqn:unitary}).
Conversely,
(\ref{eqn:block}) defines the unitary operator $U$ from ${\cal H}_{B}$ to ${\cal H}_{F}$
without the use of representations of Cuntz algebras.
%(iii)
\item
From Theorem \ref{Thm:main}(ii),
number operators $N_{B}$ and $N_{F}$ of bosons and fermions 
on their Fock spaces are
transformed as $UN_{B}=N_{F}U$ \cite{BR}.
%(iv)
\item
Theorem \ref{Thm:main} is an example of
branching law (without nontrivial branch) of $\co{2}$, ${\cal A}$ and ${\cal B}$.
We will show other branching laws in $\S$ \ref{subsection:thirdfive}.
\end{enumerate}
}
\end{rem}

%
% Remark 1.6
%
\begin{rem}
\label{rem:six}
{\rm
Since every dimensions of representation spaces 
${\cal H}_{B},{\cal H}_{F},{\cal H}$ in $\S$ \ref{subsection:firstthree}
are countably infinite,
it is clear that
there exists a one-to-one correspondence among their
standard orthonormal basis. 
Hence there exist many choices of the correspondence and 
there is no criterion to choose a correspondence in general.
On the other hand, 
the formula (\ref{eqn:block}) is derived from
two branching laws of representations associated with the RBS and the RFS.
It is surprising that pure representation theoretical 
results give unique correspondence of state vectors with 
the physical statement about the particle number conservation.
}
\end{rem}

In $\S$ \ref{section:second},
we will introduce representations of $\coni,\co{2},
{\cal B}$ and ${\cal A}$.
In $\S$ \ref{section:third},
we will explain branching laws 
and show relations among representations in $\S$ \ref{section:second}.
In $\S$ \ref{subsection:thirdfive},
we will show the proof of Theorem \ref{Thm:main}.
In $\S$ \ref{section:fourth},
we will show examples.

%%%%%%%%%%%%%%%%%%%%%%%%%%%%%%%%%%%%%%%%%%%%%%%%%%%%%%%%%%%%%%%%%%%
% 
% Section 2
%
\sftt{Representations of algebras}
\label{section:second}
In this section, we introduce representations of $\coni$, $\co{2}$, ${\cal B}$
and ${\cal A}$ independently and show their properties.
%%%%%%%%%%%%%%%%%%%%%%%%%%%%%%%%%%%%%%%%%%%%%%%%%%%%%%%%%%%%%%
%
% subsection 2.1
%
\ssft{Permutative representations of Cuntz algebras}
\label{subsection:secondone}
For $N=2,3,\ldots,+\infty$, 
let $\con$ denote the {\it Cuntz algebra} \cite{C}, that is, a C$^{*}$-algebra 
which is universally generated by $s_{1},\ldots,s_{N}$ satisfying
$s_{i}^{*}s_{j}=\delta_{ij}I$ for $i,j\edot$ and
\[\sum_{i=1}^{N}s_{i}s_{i}^{*}=I\quad(\mbox{if } N<+\infty),\quad
\sum_{i=1}^{k}s_{i}s_{i}^{*}\leq I,\quad k= 1,2,\ldots\quad
(\mbox{if }N = +\infty)\]
where $I$ denotes the unit of $\con$.

In $\S$ 2.2 of \cite{GG},
a Cuntz algebra-like object appears.
Gopakumar and Gross call it the Cuntz algebra. 
They regard that this algebra corresponds with (Maxwell-) Boltzmann statistics
(\cite{Sakurai}, p 362).
According to their interpretation, we illustrate relations
between algebras and statistics as follows:\\

\begin{tabular}{c|c|c|c}
\hline
\textsf{algebra} & boson & fermion & $\con$ \\
\hline
\textsf{statistics} & Bose-Einstein & Fermi-Dirac & Maxwell-Boltzmann(?)\\
\hline
\end{tabular}
\\

Since $\con$ is simple, that is, there is no
nontrivial closed two-sided ideal,
any unital homomorphism from $\con$ to a C$^{*}$-algebra is injective.
If $t_{1},\ldots,t_{n}$ are elements of a unital C$^{*}$-algebra
${\goth A}$ such that
$t_{1},\ldots,t_{n}$ satisfy the relations of canonical generators of $\con$,
then the correspondence $s_{i}\mapsto t_{i}$ for $i\edot$
is uniquely extended to a $*$-embedding
of $\con$ into ${\goth A}$ from the uniqueness of $\con$.
Therefore we call such a correspondence 
among generators by an embedding of $\con$ into ${\goth A}$.

Define $X_{N}\equiv \{1,\ldots,N\}$ for $2\leq N<\infty$
and $X_{\infty}\equiv {\bf N}$.
For $N=2,\ldots,\infty$ and $k=1,\ldots,\infty$, 
define the product set $X_{N}^{k}\equiv (X_{N})^{k}$ of $X_{N}$.
Let $\{s_{n}:n\in X_{N}\}$ denote the set of canonical generators of $\con$
for $2\leq N\leq \infty$.
%
% Definition 2.1
% 
\begin{defi}
\label{defi:first}
\begin{enumerate}
%(i)
\item
A representation $({\cal H},\pi)$ of $\con$
is permutative if 
there exists an orthonormal basis ${\cal E}\equiv \{e_{n}:n\in \Lambda\}$ 
of ${\cal H}$
such that $\pi(s_{i}){\cal E}\subset {\cal E}$ for each $i\in X_{N}$
\cite{BJ,DaPi2,DaPi3}.
%(ii)
\item
For $J=(j_{l})_{l=1}^{k}\in X_{N}^{k}$ with $1\leq k < \infty$,
let $P_{N}(J)$ denote the class of representations $({\cal H}, \pi)$ of $\con$ 
with a cyclic unit vector $\Omega\in {\cal H}$
such that $\pi(s_{J})\Omega=\Omega$
and $\{\pi(s_{j_{l}}\cdots s_{j_{k}})\Omega\}_{l=1}^{k}$
is an orthonormal family in ${\cal H}$
where $s_{J}\equiv s_{j_{1}}\cdots s_{j_{k}}$.
\end{enumerate}
\end{defi}

\noindent
We call the vector $\Omega$ in Definition \ref{defi:first}
by the {\it GP vector} of $({\cal H},\pi)$.

Results of these classes are shown as follows.
For any $J$, $P_{N}(J)$ 
is a class of permutative representations, which 
contains only one unitary equivalence class.
From this, we can always identify $P_{N}(J)$ with a representative of $P_{N}(J)$.
The class $P_{N}(J)$ is equivalent to $P_{N}(\sigma J)$
where $\sigma J=(j_{\sigma(1)},\ldots,j_{\sigma(k)})$ for any
cyclic permutation $\sigma\in {\bf Z}_{k}$ for $J=(j_{1},\ldots,j_{k})$.
The class $P_{N}(J)$ is irreducible if and only if $\sigma J\ne J$ for any cyclic permutation
$\sigma\ne id$ \cite{BJ,DaPi2,DaPi3,K1}.
Let $P_{N}(j,\ldots,j,k)$ ($p$-times $j$, and $k$) denote $P_{N}(j^{p}k)$ for description of 
simplicity.
We summarize our results as follows.
%
% Lemma 2.2
%
\begin{lem}
\label{lem:permutativetwo}
\cite{BJ,K1}
Let ${\cal R}\equiv \{P_{2}(2^{p-1}1), P_{2}(1^{q-1}2):p,q\geq 1,\,q\ne 2\}$.
\begin{enumerate}
%(i)
\item
Any two of classes in ${\cal R}$ are not unitarily equivalent.
%(ii)
\item
All of classes in ${\cal R}$ are irreducible.
\end{enumerate}
\end{lem}

%
% Lemma 2.3
%
\begin{lem}
(\cite{RBS01}, Lemma 2.2)
\label{lem:permutative}
Let ${\cal T}\equiv \{P_{\infty}(p),P_{\infty}(1^{p}2):p\geq 1\}$.
\begin{enumerate}
%(i)
\item
Any two of classes in ${\cal T}$ are not unitarily equivalent.
%(ii)
\item
All of classes in ${\cal T}$ are irreducible.
\end{enumerate}
\end{lem}

%%%%%%%%%%%%%%%%%%%%%%%%%%%%%%%%%%%%%%%%%%%%%%%%%%%%%
%
% subsection 2.2
%
\ssft{Representations of bosons}
\label{subsection:secondtwo}
We summarize several representations of bosons and their properties.
Let ${\cal B}$ denote the $*$-algebra generated by $\{b_{n}:n\in {\bf N}\}$
which satisfies (\ref{eqn:boson}).
The algebra ${\cal B}$ is called the {\it Heisenberg algebra} \cite{Kac,MJD},
the universal enveloping algebra of the Heisenberg Lie algebra \cite{Kac},
or the {\it Weyl algebra} \cite{KP}.
A {\it representation} of ${\cal B}$
is a pair $({\cal H},\pi)$ such that ${\cal H}$ is a complex Hilbert space 
with a dense subspace ${\cal D}$ and 
$\pi$ is a $*$-homomorphism from ${\cal B}$ to the $*$-algebra 
$\{x\in {\rm End}_{{\bf C}}({\cal D}): x^{*}{\cal D}\subset{\cal D}\}$.
A {\it cyclic vector} of $({\cal H},\pi)$ is a vector $\Omega\in {\cal D}$
such that $\pi({\cal B})\Omega={\cal D}$.
%
% Definition 2.4
%
\begin{defi}
\label{defi:bosonrep}
For $\lambda,q\geq 1$ and $i=1,\ldots,q$,
let $BF_{q,i}(\lambda)$ denote
the class of representations $({\cal H},\pi)$ of ${\cal B}$
with a cyclic vector $\Omega$ satisfying
%
% Equation 2.1
%
\begin{equation}
\label{eqn:lambda}
\pi(b_{q(n-1)+i}b_{q(n-1)+i}^{*})\Omega=\lambda \Omega,\quad
\pi(b_{q(n-1)+j})\Omega=0
\end{equation}
for $n\in {\bf N}$ and $j=1,\ldots,q,\,j\ne i$.
\end{defi}
The classes $BF_{1,1}(p)$,
$BF_{2,1}(2)$ and $BF_{2,2}(2)$ are
same as $F_{p},F_{21}$ and $F_{12}$ in Definition 2.2 of \cite{RBS01}.
A representation $({\cal H},\pi)$ of ${\cal B}$ is called {\it irreducible}
if there exists a $\pi({\cal B})$-invariant dense subspace ${\cal D}$ of ${\cal H}$
such that 
if a linear operator $y$ from ${\cal D}$ to ${\cal D}$ satisfying
$y\pi(x)=\pi(x)y$ on ${\cal D}$ for any $x\in {\cal B}$,
is a scalar multiples of $I$.

%
% Lemma 2.5
%
\begin{lem}(\cite{RBS01}, Lemma 2.2)
\label{lem:bosonresult}
Let ${\cal S}\equiv \{BF_{1,1}(p),BF_{1,2}(2),BF_{2,1}(2):p\geq 1\}$.
\begin{enumerate}
%(i)
\item
For each $S\in {\cal S}$,
any two representations belonging to $S$ are unitarily equivalent.
From this, we can identify a representation belonging to $S\in {\cal S}$ 
with $S$.
%(ii)
\item
Any two of classes in ${\cal S}$
are not unitarily equivalent.
%(iii)
\item
All of classes in ${\cal S}$ are irreducible.
\end{enumerate}
\end{lem}
We see that $BF_{1,1}(1)$ is 
the Fock representation of ${\cal B}$ with the vacuum $\Omega$.

%%%%%%%%%%%%%%%%%%%%%%%%%%%%%%%%%%%%%%%%%%%%%%%%
%
% subsection 2.3
%
\ssft{Representations of fermions}
\label{subsection:secondthree}
Let ${\cal A}$ denote the $*$-algebra generated by $\{a_{n}:n\in {\bf N}\}$
which satisfies (\ref{eqn:fermion}).
The algebra ${\cal A}$ 
is isomorphic to the {\it Clifford algebra} \cite{Kac,MJD}.
The C$^{*}$-algebra ${\sf A}$ 
universally generated by ${\cal A}$ is called the {\it CAR algebra} \cite{BR}. 
Every $*$-representation ${\cal A}$
is uniquely extended to the unital $*$-representation of ${\sf A}$.
The C$^{*}$-algebra ${\sf A}$ is simple.
%
% Definition 2.6
%
\begin{defi}
\label{defi:second}
\begin{enumerate}
%(i)
\item
For $p\geq 1$ and $1\leq i\leq p$, 
let $FF_{p,i}$ denote the class of representations $({\cal H},\pi)$ of ${\cal A}$ 
with a cyclic  vector $\Omega$ satisfying
%
% Equation 2.2
%
\begin{equation}
\label{eqn:fone}
\pi(a_{p(n-1)+i})\Omega=\pi(a_{p(n-1)+j}^{*})\Omega=0
\end{equation}
for $n\geq 1,\,j=1,\ldots,p,\,j\ne i$.
%
%(ii)
\item
For $p\geq 1$ and $1\leq i\leq p$, 
let $FF_{p,i}^{*}$ denote the class of representations $({\cal H},\pi)$ 
of ${\cal A}$ 
with a cyclic vector $\Omega$ satisfying
%
% Equation 2.3
%
\begin{equation}
\label{eqn:ftwo}
\pi(a_{p(n-1)+i}^{*})\Omega=\pi(a_{p(n-1)+j})\Omega=0\quad
\end{equation}
for $n\geq 1,\,j=1,\ldots,p,\,j\ne i$.
\end{enumerate}
\end{defi}

\noindent
By definition,
$FF_{2,1}=FF_{2,2}^{*}$ and $FF_{2,2}=FF_{2,1}^{*}$.
%
% Example 2.7
%
\begin{ex}
\label{ex:fermionrep}
{\rm
\begin{enumerate}
%(i)
\item
When $p=i=1$,
\[\pi(a_{n})\Omega=0\quad(n\geq 1).\]
Hence $FF_{1,1}$ is the Fock representation of ${\cal A}$ 
with the vacuum $\Omega$.
Every representation $FF_{p,i}$ and $FF_{p,i}^{*}$ are obtained from
the Bogoliubov transformation of the $FF_{1,1}$ \cite{BR}.
Hence $FF_{p,i}$ and $FF_{p,i}^{*}$ are often called the Fock representation
in the broad sense of the term even if $(p,i)\ne (1,1)$ \cite{MJD}.
In this paper, we call only $FF_{1,1}$ the Fermi-Fock representation 
of ${\cal A}$.
%(ii)
\item
When $p=2$, $FF_{2,1}$ satisfies
\[\pi(a_{2n-1})\Omega=\pi(a_{2n}^{*})\Omega=0\quad(n\geq 1).\]
Hence $FF_{2,1}$ is the infinite wedge representation of ${\cal A}$ 
by Proposition 3.6 in \cite{IWF01}.
In the same way,
$FF_{2,2}$ is the dual infinite wedge representation of ${\cal A}$.
\end{enumerate}
}
\end{ex}

%
% Lemma 2.8
%
\begin{lem}
\label{lem:fermirep}
Let ${\cal T}\equiv \{FF_{p,i}:p\geq 1,\,1\leq i\leq p\}$.
\begin{enumerate}
%(i)
\item
For each $T\in {\cal T}$,
any two representations belonging to $T$ are unitarily equivalent.
From this, we can identify a representation belonging to $T\in {\cal T}$ 
with $T$.
%(ii)
\item
Any two of classes in ${\cal T}$
are not unitarily equivalent.
%(iii)
\item
Any class in ${\cal T}$ is irreducible.
\end{enumerate}
\end{lem}
%
% Proof
%
\pr
(i)
From Definition 2.2 and 
Theorem 2.3 of \cite{AK02R},
the statement holds.

\noindent
(ii)
We see that ${\cal A}$ is a dense $*$-subalgebra
of the fixed-point subalgebra $\co{2}^{U(1)}$
with respect to the gauge action of $\co{2}$.
From Theorem 2.3 in \cite{AK02R}
and Proposition \ref{prop:fermioneq},
we see that $FF_{p,i}$ is $P[\sigma (2^{p-1}1)]$ for a
certain $\sigma\in {\bf Z}_{p}$.
Therefore the statement holds.

\noindent
(iii) Since any class in ${\cal T}$ is given as a Bogoliubov transformation 
from the Fock representation, and 
the Fock representation is irreducible, the statement holds.
\qedh

%%%%%%%%%%%%%%%%%%%%%%%%%%%%%%%%%%%%%%%%%%%%%%%%%%%%%%%%%%
%
% Section 3
%
\sftt{Branching laws}
\label{section:third}
In this section, we show branching laws among representations
of $\coni$, $\co{2}$, ${\cal B}$ and ${\cal A}$ 
in $\S$ \ref{section:second}.
%%%%%%%%%%%%%%%%%%%%%%%%%%%%%%%%%%%%%%%%%%%%%%%%%%%%%%%%%%
%
% subsection 3.1
%
\ssft{Introduction to branching laws}
\label{subsection:thirdone}
First, we explain the notion of branching law.
For a group $G$, if there exists an embedding of $G$ into 
some other group $G^{'}$,
then any representation $\pi$ of $G^{'}$ induces 
the restriction $\pi|_{G}$ of $\pi$ on $G$.
The representation $\pi|_{G}$ is not irreducible in general
even if $\pi$ is irreducible.
If $\pi|_{G}$ is decomposed into the direct sum of
a family $\{\pi_{\lambda}:\lambda\in\Lambda\}$ of irreducible representations 
of $G$, then the equation
%
% Equation 3.1
%
\begin{equation}
\label{eqn:group}
\pi|_{G}=\bigoplus_{\lambda\in\Lambda}\pi_{\lambda}
\end{equation}
is called the {\it branching law} of $\pi$.
Especially, if $\pi|_{G}$ itself is irreducible,
then the branching law of $\pi|_{G}$ has no nontrivial branch.
The branching law can be also considered for a pair $(A,B)$
of a subalgebra $A$ and an algebra $B$.
We can consider branching laws for the following pairs:
\begin{enumerate}
%(i)
\item
$(\coni,\co{2})$, 
%(ii)
\item
$({\cal A},\co{2})$,
%(iii)
\item
$({\cal B},\coni)$,
%(iv)
\item
$({\cal B},\co{2})$
\end{enumerate}
where ${\cal B}$ is {\it neither} a subalgebra of $\coni$
{\it nor} that of $\co{2}$, but
branching laws of permutative representations can be considered 
as if ${\cal B}$ was a subalgebra of $\coni$ by \cite{RBS01},
and we write this inclusion like relation as 
the symbol ``$\rightsquigarrow$" in the left part of the following diagram.
From $(\coni,\co{2})$ and $({\cal B},\coni)$,
we can consider branching laws for the pair $({\cal B},\co{2})$.
We illustrate relations among them as follows:

%%%%%%%%%%%%%%%%%%%%%%%%%%%%%%%%%%%%%%%%
\def\diagrams{
\put(200,0){$\co{2}$,}
\put(0,0){$\coni$}
\put(10,130){${\cal B}$}
\put(10,-110){${\cal A}$}
\put(100,0){$\hookrightarrow$}
%\put(10,100){\rotatebox{-90}{$\hookrightarrow$}}
%\put(10,100){\rotatebox{-90}{$\rightarrowtail$}}
\put(10,100){\rotatebox{-90}{$\rightsquigarrow$}}
\put(90,80){\rotatebox{-30}{$\rightsquigarrow$}}
\put(90,-80){\rotatebox{30}{$\hookrightarrow$}}
\put(60,40){$\circlearrowleft$}
}

\def\restriction{
\put(0,0){${\rm Rep}\co{2}\ni \pi$}
\put(300,0){$\pi|_{\coni}\in{\rm Rep}\coni$}
\put(300,-110){$\pi|_{{\cal A}}\quad \in{\rm Rep}{\cal A}$}
\put(300,130){$\pi|_{{\cal B}}\quad \in{\rm Rep}{\cal B}$}
\put(220,0){$\longmapsto$}
\put(320,50){\rotatebox{90}{$\rightarrowtail$}}
\put(220,-50){\rotatebox{-30}{$\longmapsto$}}
\put(220,70){\rotatebox{30}{$\rightarrowtail$}}
\put(240,40){$\circlearrowleft$}
}

\thicklines
%\framebox{
\setlength{\unitlength}{.1mm}
\begin{picture}(1000,300)
\put(100,120){\diagrams}
\put(500,120){\restriction}
\end{picture}
%}

\noindent
In consequence, we can compare two branching laws for a permutative representation
of $\co{2}$ associated with  $({\cal B},\co{2})$ 
and $({\cal A},\co{2})$.

%%%%%%%%%%%%%%%%%%%%%%%%%%%%%%%%%%%%%%%%%%%%%%%%%%%%%%%%%%
%
% subsection 3.2
%
\ssft{$\coni$ and $\co{2}$}
\label{subsection:thirdtwo}
Let $\{s_{n}:n\in {\bf N}\}$ and $\{t_{1},t_{2}\}$ 
be as in (\ref{eqn:oinfty}) and (\ref{eqn:otwo}).
Assume that $\coni$ is embedded into $\co{2}$ by (\ref{eqn:embedding}).
From this,
$t_{2}s_{n}=s_{n+1}$ for each $n\in {\bf N}$.
%
% Proposition 3.1
%
\begin{prop}
\label{prop:oni}
For $P_{N}(J)$ in Definition \ref{defi:first}(ii),
let $P_{2}(J)|_{\coni}$ denote the 
restriction of the (class of) representation $P_{2}(J)$ on $\coni$.
\begin{enumerate}
%(i)
\item
For $p\geq 1$,
%
% Equation 3.3
%
\begin{equation}
\label{eqn:threeonetwo}
P_{2}(12^{p-1})|_{\coni}=P_{\infty}(p)
\end{equation}
where we define $12^{0}=1$ for convenience.
%(ii)
\item
For $q\geq 1$, 
%
% Equation 3.4
%
\begin{equation}
\label{eqn:threeonethree}
P_{2}(1^{q}2)|_{\coni}=P_{\infty}(1^{q-1}2)
\end{equation}
where we define $1^{0}2=2$ for convenience.
\end{enumerate}
\end{prop}
%
% Proof
%
\pr
(i)
Fix $p\geq 1$.
Let $({\cal H},\pi)$ be $P_{2}(12^{p-1})$ with the GP vector $\Omega$.
We identify $\pi(t_{i})$ and $t_{i}$ for $i=1,2$.
Then $t_{1}t_{2}^{p-1}\Omega=\Omega$.
Define $\Omega^{'}\equiv t_{2}^{p-1}\Omega$.
Then $\Omega^{'}$ is a cyclic vector and $s_{p}\Omega^{'}=\Omega^{'}$.
Let $V\equiv {\rm Lin}\langle\{t_{J}\Omega^{'}:J\in\{1,2\}^{*}\}\rangle$
where $\{1,2\}^{*}\equiv \bigcup_{l\geq 1}\{1,2\}^{l}$.
Then we see that $V$ is a dense subspace of ${\cal H}$.

Fix $p\geq 1$.
For $n_{2},\ldots,n_{k},m_{1},\ldots,m_{k-1}\geq 1$ and $n_{1},m_{k}\geq 0$,
let $J=(1^{n_{1}}2^{m_{1}}\cdots 1^{n_{k}}2^{m_{k}})\in \{1,2\}^{*}$.
Then
%
% Equation 3.5
%
\begin{equation}
\label{eqn:threeonefour}
t_{J}\Omega^{'}=s_{1}^{n_{1}}\cdot s_{m_{1}+1}s_{1}^{n_{2}-1}\cdot
s_{m_{2}+1}s_{1}^{n_{3}-1}
\cdots 
s_{m_{k-1}+1}s_{1}^{n_{k}-1}s_{m_{k}+p}\Omega^{'}.
\end{equation}
From this,
$t_{J}\Omega^{'}\in \coni\Omega^{'}$ for each $J\in\{1,2\}^{*}$.
Therefore $V\subset \coni\Omega^{'}$.
This implies ${\cal H}=\overline{\coni\Omega^{'}}$.
Hence $({\cal H},\pi|_{\coni})$ is $P_{\infty}(p)$.

\noindent
(ii)
Assume that a vector $\Omega$ satisfies $t_{1^{q}2}\Omega^{'}=\Omega^{'}$.
Define $\Omega\equiv t_{1^{q-1}2}\Omega^{'}$.
Then
$s_{1}^{q-1}s_{2}\Omega=t_{1}^{q-1}t_{2}t_{1}\Omega=\Omega$.
Hence $V\equiv \overline{\coni\Omega}$ is $P_{\infty}(1^{q-1}2)$.
On the other hand,
$t_{1}\Omega=s_{1}\Omega$ and 
$t_{2}\Omega=
t_{2}t_{1^{q-1}2}t_{1}\Omega=s_{2}s_{1^{q-2}}s_{2}\Omega$.
From these, $t_{i}\Omega\in V$ for $i=1,2$.
Furthermore, we see that
%
% Equation 3.6
%
\begin{equation}
\label{eqn:threeonesix}
t_{1}s_{j_{1}}\cdots s_{j_{k}}\Omega=
s_{1}s_{j_{1}}\cdots s_{j_{k}}\Omega,\quad
t_{2}s_{j_{1}}\cdots s_{j_{k}}\Omega=
s_{j_{1}+1}s_{j_{2}}\cdots s_{j_{k}}\Omega.
\end{equation}
Therefore $t_{J}\Omega\in V$ for each $J\in\{1,2\}^{*}$.
Since $\Omega$ is a cyclic vector of $P_{2}(1^{q}2)$,
the statement holds.
\qedh

%%%%%%%%%%%%%%%%%%%%%%%%%%%%%%%%%%%%%%%%%%%%%%%%%%%%%%%%%%%%
%
% subsection 3.3
%
\ssft{Bosons and $\co{2}$}
\label{subsection:thirdthree}
We show relations among representations of boson and $\co{2}$ 
according to (\ref{eqn:rbszero}) and (\ref{eqn:embedding}).
From (\ref{eqn:rbszero}), 
%
% Equation 3.6
%
\begin{equation}
\label{eqn:commboson}
s_{n}b_{m}=b_{m+1}s_{n},\quad s_{n}b_{m}^{*}=b_{m+1}^{*}s_{n}
\quad (n,m\in {\bf N}).
\end{equation}
%
%
% Proposition 3.2
%
\begin{prop}
\label{prop:bosoneq}
\begin{enumerate}
%(i)
\item
For $p\geq 1$, 
$P_{2}(12^{p-1})|_{{\cal B}}=BF_{1,1}(p)$.
%(ii)
\item
For $p\geq 1$,	
$P_{\infty}(1^{p-1}2)|_{{\cal B}}=BF_{p,1}(2)\oplus \cdots \oplus BF_{p,p}(2)$.
%(iii)
\item
For $p\geq 1$,	
$P_{2}(1^{p}2)|_{{\cal B}}=BF_{p,1}(2)\oplus \cdots \oplus BF_{p,p}(2)$.
\end{enumerate}
\end{prop}
%
% Proof
%
\pr
(i) From Proposition \ref{prop:oni}(i)
and $P_{\infty}(p)|_{{\cal B}}=BF_{1,1}(p)$ 
by Theorem 1.1 of \cite{RBS01}, 
the statement holds.

\noindent
(ii)
When $p=1$, the statement is proved by Theorem 1.1 of \cite{RBS01}.
Assume $p\geq 2$.
Let $\Omega$ be a non-zero vector satisfying $s_{1}^{p-1}s_{2}\Omega=\Omega$.
Define
%
% Equation 3.9
%
\begin{equation}
\label{eqn:threenine}
\Omega_{1}\equiv \Omega,\quad
\Omega_{2}\equiv s_{1}^{p-2}s_{2}\Omega,\ldots,\quad
\Omega_{p-1}\equiv s_{1}s_{2}\Omega,\quad \Omega_{p}\equiv s_{2}\Omega
\end{equation}
and $T_{i}\equiv s_{1}^{p-i}s_{2}s_{1}^{i-1}$ where $s_{1}^{0}\equiv I$.
Then $T_{i}\Omega_{i}=\Omega_{i}$ for each $i=1,\ldots,p$.
We see that
$b_{p-j+1}\Omega_{i}=\delta_{i,j}s_{1}^{p}\Omega_{i}$ for $i,j=1,\ldots,p$.
From this,
%
% Equation 3.9
%
\begin{equation}
\label{eqn:threeeleven}
b_{p(n-1)+p-j+1}\Omega_{i}=\delta_{ij}
T_{i}^{n-1}s_{1}^{p}\Omega_{i}\quad (i,j=1,\ldots,p,\,n\geq 1).
\end{equation}
Moreover
$b_{p(n-1)+p-i+1}^{*}b_{p(n-1)+p-i+1}\Omega_{i}=\Omega_{i}$
for $i=1,\ldots,p,\,n\geq 1$.
Hence $V_{i}\equiv {\cal B}\Omega_{i}$ is $BF_{p,p-i+1}(2)$ 
for each $i=1,\ldots,p$.
Therefore
$W\equiv V_{1}\oplus\cdots\oplus V_{p}$ is 
$BF_{p,1}(2)\oplus \cdots \oplus BF_{p,p}(2)$.

On the other hand,
%
% Equation 3.10
%
\begin{equation}
\label{eqn:threethirteen}
s_{1}\Omega_{1}=b_{1}\Omega_{p},\quad
s_{1}\Omega_{i}=\Omega_{i-1},\quad s_{2}\Omega_{1}=\Omega_{p},\quad
s_{2}\Omega_{i}=b_{1}^{*}\Omega_{i-1},
\end{equation}
%
% Equation 3.11
%
\begin{equation}
\label{eqn:threefourteen}
s_{n}\Omega_{1}=\{(n-2)!\}^{-1/2}(b_{1}^{*})^{n-2}\Omega_{p},\quad
s_{n}\Omega_{i}=\{(n-1)!\}^{-1/2}(b_{1}^{*})^{n-1}\Omega_{i-1}
\end{equation}
for $n\geq 3$ and $i=2,\ldots,p$.
These imply that 
$s_{n}\Omega_{i}\in W$ for each $i=1,\ldots,p$ and $n\in {\bf N}$.
From this and (\ref{eqn:commboson}),
$s_{n}(a_{n_{1}}^{*})^{k_{1}}\cdots (a_{n_{r}}^{*})^{k_{r}}
a_{m_{1}}^{l_{1}}\cdots a_{m_{t}}^{l_{t}}\Omega_{i}\in W$
for each $i=1,\ldots,p$, $n_{1},\ldots,n_{r}$, 
$m_{1},\ldots,m_{t}$,
$k_{1},\ldots,k_{r}$ and $l_{1},\ldots,l_{t}$.
Hence 
$s_{n}W\subset W$ for each $n\in {\bf N}$.
From this, $s_{J}\Omega=s_{J}\Omega_{1}\in W$ for each $J$.
Therefore $W$ is dense in $\coni\Omega$.
This implies the statement.

\noindent
(iii)
From (ii) and Proposition \ref{prop:oni}(ii), the statement holds.
\qedh

%%%%%%%%%%%%%%%%%%%%%%%%%%%%%%%%%%%%%%%%%%%%%%%%%%%%
%
% subsection 3.4
%
\ssft{Fermions and $\co{2}$}
\label{subsection:thirdfour}
We show relations among representations of fermions and $\co{2}$.
From (\ref{eqn:rfszero}), 
%
% Equation 3.11
%
\begin{equation}
\label{eqn:commfermion}
t_{i}a_{m}=(-1)^{i-1}a_{m+1}t_{i},\quad (-1)^{i-1}t_{i}a_{m}^{*}=a_{m+1}^{*}t_{i}
\end{equation}
for $i=1,2$ and $m\in {\bf N}$.
For $p\geq 1$,
let $({\cal H},\pi)$ be $P_{2}(2^{p-1}1)$ with the GP vector $\Omega$.
Define
%
% Equation 3.12
%
\begin{equation}
\label{eqn:omega}
\Omega_{1}\equiv \Omega,\quad \Omega_{2}\equiv t_{2}^{p-2}t_{1}\Omega,
\ldots,\quad \Omega_{p-1}\equiv t_{2}t_{1}\Omega,\quad 
\Omega_{p}\equiv t_{1}\Omega.
\end{equation}
From (\ref{eqn:commfermion}), the following is verified
%
% Lemma 3.3
%
\begin{lem}
\label{lem:fermionone}
\begin{enumerate}
%(i)
\item
When $p=1$,
$a_{n}^{*}\Omega_{1}=t_{1}^{n-1}t_{2}\Omega_{1}$
for each $n\geq 1$.
%(ii)
\item
When $p\geq 2$,
$a_{p-i+1}^{*}\Omega_{j}=\delta_{ij}(-1)^{p-i}t_{2}^{p-i+1}\Omega$
for $i,j=1,\ldots,p$.
%(iii)
\item
When $p\geq 2$,
for $i,j=1,\ldots,p,\,l\geq 2$,
%
% Equation 3.13
%
\begin{equation}
\label{eqn:threetwenty}
a_{p(l-1)+p-i+1}^{*}\Omega_{j}
=\delta_{ij}(-1)^{p-i+(p-1)(l-1)}t_{2}^{p-i}t_{1}(t_{2}^{p-1}t_{1})^{l-2}
t_{2}^{p}\Omega.
\end{equation}
In particular, when $i=j=1$,
%
% Equation 3.14
%
\begin{equation}
\label{eqn:threetwentyone}
a_{pl}^{*}\Omega= (-1)^{(p-1)l}
(t_{2}^{p-1}t_{1})^{l-1}t_{2}^{p}\Omega,\quad 
a_{p(l-1)+i}^{*}\Omega=0
\end{equation}
for $l\geq 1,\,1\leq i\leq p-1$.
\end{enumerate}
\end{lem}

\noindent
From Lemma \ref{lem:fermionone} and (\ref{eqn:fermion}),
the following holds.
%
% Lemma 3.4
%
\begin{lem}
\label{lem:fermiontwo}
\begin{enumerate}
%(i)
\item
For $i,j=1,\ldots,p,\,i\ne j,\,l\geq 1$,
%
% Equation 3.15, 3.16
%
\begin{eqnarray}
\label{eqn:threetwentytwo}
a_{p(l-1)+p-i+1}a_{p(l-1)+p-i+1}^{*}\Omega_{i}=\Omega_{i},\\
\nonumber
\\
\label{eqn:threeb}
a_{p(l-1)+p-i+1}^{*}a_{p(l-1)+p-i+1}\Omega_{j}=\Omega_{j}.
\end{eqnarray}
In particular,
for $l\geq 1$ and $1\leq i\leq p-1$,
%
% Equation 3.17
%
\begin{equation}
\label{eqn:threetwentythree}
a_{pl}a_{pl}^{*}\Omega= \Omega,\quad 
a_{p(l-1)+i}^{*}a_{p(l-1)+i}\Omega=\Omega.
\end{equation}
%
%(ii)
\item
For $i,j=1,\ldots,p,\,i\ne j,\,l\geq 1$,
%
% Equation 3.24
%
\begin{equation}
\label{eqn:threetwentyfour}
a_{p(l-1)+p-i+1}\Omega_{i}=a_{p(l-1)+p-i+1}^{*}\Omega_{j}=0.
\end{equation}
\end{enumerate}
\end{lem}

%
% Lemma 3.6
%
\begin{lem}
\label{lem:fermiaction}
For $\{\Omega_{j}\}_{j=1}^{p}$ in (\ref{eqn:omega}), the following holds:
\[t_{1}\Omega_{1}=\Omega_{p},\quad t_{2}\Omega_{1}=a_{1}^{*}\Omega_{p},\quad
t_{1}\Omega_{j}=a_{1}\Omega_{j-1},\quad
t_{2}\Omega_{j}=\Omega_{j-1}\quad(2\leq j\leq p).\]
\end{lem}
%
% Proof
%
\pr
By definition, $t_{1}\Omega_{1}=\Omega_{p}$ and 
$t_{2}\Omega_{j}=\Omega_{j-1}$ for $2\leq j\leq p$.
On the other hand,
$t_{2}\Omega_{1}=t_{2}t_{1}^{*}\Omega_{p}=a_{1}^{*}\Omega_{p}$ and 
$t_{1}\Omega_{j}=t_{1}t_{2}^{*}\Omega_{j-1}=a_{1}\Omega_{j-1}$ 
for $2\leq j\leq p$.
\qedh

%
% Proposition 3.8
%
\begin{prop}
\label{prop:fermioneq}
For any $p\geq 1$, the following holds:
\[P_{2}(2^{p-1}1)|_{{\cal A}}=FF_{p,1}\oplus\cdots \oplus FF_{p,p},\quad
P_{2}(1^{p-1}2)|_{{\cal A}}=FF_{p,1}^{*}\oplus\cdots \oplus FF_{p,p}^{*}.\]
\end{prop}
%
% Proof
%
\pr
Let ${\cal H}$ be the representation space of $P_{2}(2^{p-1}1)$.
By Lemma \ref{lem:fermionone},
we see that 
$\overline{{\cal A}\Omega_{j}}$ is $FF_{p,j}$
for each $j=1,\ldots,p$.
Therefore
$FF_{p,1}\oplus\cdots \oplus FF_{p,p}$ is a subrepresentation of 
${\cal H}$.
It is sufficient to show that 
$W\equiv {\cal A}\Omega_{1}\oplus\cdots\oplus {\cal A}\Omega_{p}$
is dense in ${\cal H}$.

By Lemma \ref{lem:fermiaction},
$t_{i}\Omega_{j}\in W$ for each $i=1,2$ and $j=1,\ldots,p$.
From this and (\ref{eqn:commfermion}),
%
% Equation 3.25
%
\begin{equation}
\label{eqn:threetwentyfive}
t_{i}a_{S}^{*}a_{T}\Omega_{j}
=(-1)^{(|S|+|T|)(i-1)}a_{S+1}^{*}a_{T+1}t_{i}\Omega_{j}\in W
\end{equation}
where $a_{S}^{*}=a_{n_{1}}^{*}\cdots a_{n_{k}}^{*}$,
$a_{T}\equiv (a_{T}^{*})^{*}$
and $|S|=k$, $S+1\equiv \{n_{1}+1,\ldots,n_{k}+1\}$
for $S=\{n_{1},\ldots,n_{k}\}$.
Let ${\cal R}\equiv \{t_{J}\Omega:J\in \{1,2\}^{*}\}$.
From Lemma \ref{lem:fermiaction},
${\cal R}\subset W$.
Since the linear hull of ${\cal R}$ is dense
in ${\cal H}$, $W$ is also dense in ${\cal H}$.
Hence the first statement holds.

Define the automorphism $\alpha$ of $\co{2}$ 
by $\alpha(t_{1})\equiv t_{2}$ and $\alpha(t_{2})\equiv t_{1}$.
Then $P_{2}(2^{p-1}1)\circ \alpha=P_{2}(1^{p-1}2)$.
Furthermore $\alpha|_{{\cal A}}$ is also an automorphism of ${\cal A}$.
We see that $FF_{p,i}\circ \alpha=FF_{p,i}^{*}$.
From this and the first statement, the second statement holds.
\qedh

For other results of representations of fermions and Cuntz algebras,
see \cite{AK1,AK06,AK4,AK05,AK02R,IWF01}.

%%%%%%%%%%%%%%%%%%%%%%%%%%%%%%%%%%%%%%%%%%%%%%%%%%%%%%%%%%
%
% subsection 3.5
%
\ssft{Bosons and fermions}
\label{subsection:thirdfive}
In this subsection, we summarize relations among 
representations of algebras $\co{2},{\cal B}$ and ${\cal A}$,
and prove Theorem \ref{Thm:main}.
From Proposition \ref{prop:bosoneq} and \ref{prop:fermioneq},
we obtain the unitary $U_{p}$ from $BF_{1,1}(p)$ to
$FF_{p,1}\oplus \cdots \oplus FF_{p,p}$ 
for $p\geq 1$ as follows:

\noindent
%
% Equation 3.20
%
\begin{equation}
\label{eqn:up}
\thicklines
\setlength{\unitlength}{.1mm}
\begin{picture}(1000,250)(120,210)
\put(300,400){$P_{2}(2^{p-1}1)$}
\put(50,250){$BF_{1,1}(p)$}
\put(600,250){$FF_{p,1}\oplus \cdots \oplus FF_{p,p}$,}
\put(270,370){\vector(-2,-1){120}}
\put(480,370){\vector(2,-1){120}}
\put(230,270){\vector(1,0){300}}
%\put(230,270){\vector(-1,0){0}}
\put(120,360){$V_{B,p}$}
\put(570,360){$V_{F,p}$}
\put(370,290){$U_{p}$}
\put(360,340){$\circlearrowleft$}
\put(930,350){$U_{p}\equiv V_{F,p}(V_{B,p})^{*}$}
\end{picture}
\end{equation}
where $V_{B,p}$ denotes the unitary 
from $P_{2}(2^{p-1}1)|_{{\cal B}}$ to $BF_{1,1}(p)$ 
and $V_{F,p}$ denote the unitary from $P_{2}(2^{p-1}1)|_{{\cal A}}$
to $FF_{p,1}\oplus \cdots \oplus FF_{p,p}$ with respect to
their unitary equivalences.

%
% Example 3.7
%
\begin{ex}
\label{ex:eight}
{\rm
\begin{enumerate}
%(i)
\item
For Example \ref{ex:fermionrep}(i),
if $p=1$, then we obtain the unitary $U=U_{1}$ from
the Bose-Fock space to the Fermi-Fock space.
%(ii)
\item
From Example \ref{ex:fermionrep}(ii),
$U_{2}$ is the following unitary: 
\[U_{2}:BF_{1,1}(2)\to IW\oplus IW^{*}\]
where $IW$ and $IW^{*}$ denote
the infinite wedge representation 
and 
the dual infinite wedge representation of ${\cal A}$, respectively.
\end{enumerate}
}
\end{ex}

Let $t_{1},t_{2}$ denote canonical generators of $\co{2}$
and let $\Omega$ be the GP vector of $P_{2}(1)$.
For $1\leq n_{1}<\cdots<n_{m}$ and
$k_{1},\ldots,k_{m}\in {\bf N}$,
the RBS $\{b_{n}:n\in {\bf N}\}$ 
on $P_{2}(1)$ satisfies the following
Example 3.2 in \cite{RBS01}:
%
% Equation 3.26
% 
\begin{equation}
\label{eqn:rbs}
(b_{n_{1}}^{*})^{k_{1}}\cdots (b_{n_{m}}^{*})^{k_{m}}\Omega
=\prod_{i=1}^{m}\sqrt{k_{i}!}\,\,
t_{1}^{n_{1}-1}t^{k_{1}}_{2}t_{1}^{n_{2}-n_{1}}t^{k_{2}}_{2}
\cdots t_{1}^{n_{m}-n_{m-1}}t^{k_{m}}_{2}\Omega.
\end{equation}
For $1\leq k_{1}<k_{2}<\cdots <k_{l}$,
the restriction of $P_{2}(1)$ on
the RFS $\{a_{n}:n\in {\bf N}\}$ satisfies the following
(\cite{AK1}, (3.40)):
%
% Equation 3.22
% 
\begin{equation}
\label{eqn:rfs}
a_{k_{1}}^{*}\cdots a_{k_{l}}^{*}\Omega
=t_{1}^{k_{1}-1}t_{2}t_{1}^{k_{2}-k_{1}-1}t_{2}\cdots 
t_{1}^{k_{l}-k_{l-1}-1}t_{2}\Omega.
\end{equation}
From these two formulae,
it is shown that any state vector in both Fock spaces are described
by using the canonical generators of $\co{2}$ and $\Omega$.

According to formulae in $\S$ \ref{subsection:secondone},
we show the proof of Theorem \ref{Thm:main}.
\\
\\
{\it Proof of Theorem \ref{Thm:main}.}
(i)
Since $V_{B}\Omega=\Omega_{B}$
and $V_{F}\Omega=\Omega_{F}$,
we see that $U\Omega_{B}=\Omega_{F}$.
From (\ref{eqn:commfermion}),
%
% Equation 3.23
%
\begin{equation}
\label{eqn:large}
t_{1}^{n}t_{2}^{m}=A_{n,m}t_{1}^{n+m}
\quad(n,m\geq 1).
\end{equation}
From (\ref{eqn:large}), (\ref{eqn:rbs}) and (\ref{eqn:rfs}),
we obtain the following:\\
\\
$U(b_{n_{1}}^{*})^{k_{1}}\cdots (b_{n_{m}}^{*})^{k_{m}}\Omega_{B}$
\[
\begin{array}{rl}
=&C
V_{F}t_{1}^{n_{1}-1}t^{k_{1}}_{2}t_{1}^{n_{2}-n_{1}}t^{k_{2}}_{2}
\cdots t_{1}^{n_{m}-n_{m-1}}t^{k_{m}}_{2}\Omega\\
=&C
A_{n_{1}-1,k_{1}}V_{F}t_{1}^{n_{1}+k_{1}-1}t_{1}^{n_{2}-n_{1}}t^{k_{2}}_{2}
\cdots t_{1}^{n_{m}-n_{m-1}}t^{k_{m}}_{2}\Omega\\
=&C
A_{n_{1}-1,k_{1}}
A_{n_{2}+k_{1}-1,k_{2}}
V_{F}t_{1}^{n_{2}+k_{1}+k_{2}-1}t_{1}^{n_{3}-n_{2}}t^{k_{3}}_{3}
\cdots t_{1}^{n_{m}-n_{m-1}}t^{k_{m}}_{2}\Omega\\
=&\cdots \\
=&C
A_{n_{1}-1,k_{1}}\cdots 
A_{n_{m}+k_{1}+\cdots+k_{m-1}-1,k_{m}}\Omega_{F}.\\
\end{array}
\]
Hence (\ref{eqn:block}) holds.

\noindent
(ii)
From (i), we see that the particle number of both sides of (\ref{eqn:block}) 
is $k_{1}+\cdots+k_{m}$.
Hence the statement holds.
\qedh

%%%%%%%%%%%%%%%%%%%%%%%%%%%%%%%%%%%%%%%%%%%%%%%%%%%%%%%
%
% section 4
%
\sftt{Examples}
\label{section:fourth}
%%%%%%%%%%%%%%%%%%%%%%%%%%%%%%%%%%%%%%%%%%%%%%%%%%%%%%%%%%%%%%%%
%
% subsection 4.1
%
\ssft{State vectors in Fock representations}
\label{subsection:fourthone}
We show concrete examples of Theorem \ref{Thm:main}.
Relations among $n$-particle states for $n=1,2,3$ are shown as follows.
In this subsection,
we write both $\Omega_{B}$ and $\Omega_{F}$ in Theorem \ref{Thm:main}
as the same symbol $\Omega$
and omit the unitary operator $U$.
%
% Example 4.1
%
\begin{ex}
\label{ex:concrete}
{\rm
\begin{enumerate}
%(i)
\item
For $n\geq 1$,
$b_{n}^{*}\Omega=a_{n}^{*}\Omega$.
%(ii)
\item
For $n,m\in {\bf N}$, assume  $1\leq n<m$. Then
\[b^{*}_{n}b^{*}_{m}\Omega=a_{n}^{*}a_{m+1}^{*}\Omega,\quad
(b^{*}_{n})^{2}\Omega=\sqrt{2}a_{n}^{*}a_{n+1}^{*}\Omega.\]
%(iii)
\item
For $n,m,l\in{\bf N}$,
assume $1\leq n<m<l$. Then
\[
\begin{array}{rl}
b^{*}_{n}b^{*}_{m}b_{l}^{*}\Omega&=
a_{n}^{*}a_{m+1}^{*}a_{l+2}^{*}\Omega,\\
\\
b^{*}_{n}(b^{*}_{m})^{2}\Omega&=\sqrt{2}a_{n}^{*}a_{m+1}^{*}a_{m+2}^{*}\Omega,\\
\\
(b^{*}_{n})^{2}b^{*}_{m}\Omega
&=\sqrt{2}a_{n}^{*}a_{n+1}^{*}a_{m+2}^{*}\Omega,\\
\\
(b^{*}_{n})^{3}\Omega&=\sqrt{6}a_{n}^{*}a_{n+1}^{*}a_{n+2}^{*}\Omega.
\end{array}
\]
%%%%%%%%%%%%%%%%%%%%%%%%%%%%%%%%%%%%%%%%%%%%%%%%%%%%%%%%%%%%%%%%%
\end{enumerate}
}
\end{ex}

\noindent
From Example \ref{ex:concrete}, 
we see that state vectors of boson and fermion 
are similar for few particle number case. 

Next, we show more general relations.
%
% Example 4.2
%
\begin{ex}
\label{ex:general}
{\rm
\begin{enumerate}
%(i)
\item
For each $n\geq 1$ and $m\geq 0$,
\[(b_{n}^{*})^{m+1}\Omega=\{(m+1)!\}^{1/2}a_{n}^{*}a_{n+1}^{*}\cdots a_{n+m}^{*}\Omega.\]
%(ii)
\item
For $1\leq n_{1}<n_{2}<\cdots<n_{l}$,
\[b_{n_{1}}^{*}\cdots b_{n_{l}}^{*}\Omega
=a_{n_{1}}^{*}a_{n_{2}+1}^{*}\cdots a_{n_{l}+l-1}^{*}\Omega.\]
In particular,
\[b_{n}^{*}b_{n+1}^{*}\cdots b_{n+l-1}^{*}\Omega
=a_{n}^{*}a_{n+2}^{*}\cdots a_{n+2(l-1)}^{*}\Omega
\quad(n,l\geq 1).\]
\end{enumerate}
}
\end{ex}

\noindent
From Example \ref{ex:general}(i) and Theorem \ref{Thm:main},
the multiplicity of a mode of boson is nearly associated with
the length of a block of fermions in (\ref{eqn:block}).

We illustrate state vectors of representations 
of $\co{2},\coni,{\cal A}$ and ${\cal B}$ in Theorem \ref{Thm:main}
as follows:

%%%%%%%%%%%%%%%%%%%%%%%%%%%%%%%%%%%%%%%%%%%%%%%%%%%%%%%%%%%%%%%%%%%%%%%
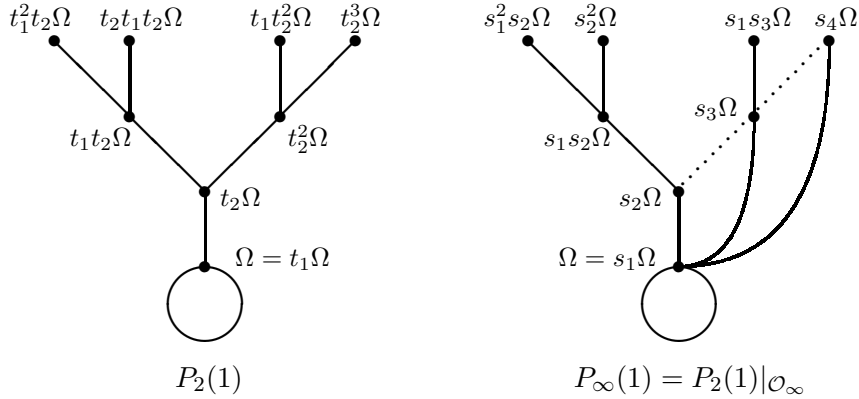
\def\vertices{
\put(0,0){$\bullet$}
\put(0,100){$\bullet$}
\put(100,200){$\bullet$}
\put(-100,200){$\bullet$}
\put(200,300){$\bullet$}
\put(100,300){$\bullet$}
\put(-200,300){$\bullet$}
\put(-100,300){$\bullet$}
}
%%%%%%%%%%%%%%%%%%%%%%%%%%%%%%%%%%%%%%%%%%%%%%%%
\def\otwo{
\put(-30,-150){$P_{2}(1)$}
\put(0,0){\vertices}
\put(10,-40){\circle{100}}
\put(10,10){\line(0,1){100}}
\put(10,110){\line(1,1){100}}
\put(10,110){\line(-1,1){100}}
\put(110,210){\line(1,1){100}}
\put(110,210){\line(0,1){100}}
\put(-90,210){\line(0,1){100}}
\put(-90,210){\line(-1,1){100}}
{\small
\put(50,10){$\Omega=t_{1}\Omega$}
\put(30,90){$t_{2}\Omega$}
\put(-170,170){$t_{1}t_{2}\Omega$}
\put(-250,330){$t_{1}^{2}t_{2}\Omega$}
\put(-130,330){$t_{2}t_{1}t_{2}\Omega$}
\put(70,330){$t_{1}t_{2}^{2}\Omega$}
\put(190,330){$t_{2}^{3}\Omega$}
\put(120,170){$t_{2}^{2}\Omega$}
}
}
%%%%%%%%%%%%%%%%%%%%%%%%%%%%%%%%%%%%%%%%%%%%%%%%
\def\fermi{
\put(-150,-60){$P_{2}(1)|_{{\cal A}}=$ Fermi-Fock}
\put(0,0){\vertices}
%\put(10,-40){\circle{100}}
\put(10,10){\line(0,1){100}}
\put(10,110){\line(1,1){100}}
%\put(10,110){\line(-1,1){100}}
\put(110,210){\line(1,1){100}}
%\put(110,210){\line(0,1){100}}
%\put(-90,210){\line(0,1){100}}
%\put(-90,210){\line(-1,1){100}}
\put(110,310){\line(1,0){100}}
\put(-90,210){\line(1,0){200}}
\qbezier(10,10)(-90,10)(-93,210)
\qbezier(10,10)(-190,10)(-190,310)
\qbezier(-190,310)(-40,410)(110,310)
\qbezier(-90,310)(60,210)(210,310)
\multiput(0,100)(-10,10){20}{$\cdot$}
\multiput(104,200)(0,10){10}{$\cdot$}
\put(10,110){\line(-1,2){100}}
\put(-190,310){\line(1,0){100}}
\put(-90,210){\line(2,1){200}}
{\small
\put(40,0){$\Omega$}
\put(30,90){$a_{1}^{*}\Omega$}
\put(-160,175){$a_{2}^{*}\Omega$}
\put(-250,330){$a_{3}^{*}\Omega$}
\put(-70,310){$a_{1}^{*}a_{3}^{*}\Omega$}
\put(80,350){$a_{2}^{*}a_{3}^{*}\Omega$}
\put(190,330){$a_{1}^{*}a_{2}^{*}a_{3}^{*}\Omega$}
\put(120,170){$a_{1}^{*}a_{2}^{*}\Omega$}
}
}
%%%%%%%%%%%%%%%%%%%%%%%%%%%%%%%%%%%%%%%%%
\def\oni{
\put(-130,-150){$P_{\infty}(1)=P_{2}(1)|_{\coni}$}
\put(0,0){\vertices}
\put(10,-40){\circle{100}}
\put(10,10){\line(0,1){100}}
\multiput(0,100)(10,10){20}{$\cdot$}
\qbezier(10,10)(110,10)(110,210)
\qbezier(10,10)(210,10)(210,310)
\put(10,110){\line(-1,1){100}}
%\put(110,210){\line(1,1){100}}
\put(110,210){\line(0,1){100}}
\put(-90,210){\line(0,1){100}}
\put(-90,210){\line(-1,1){100}}
{\small
\put(-150,10){$\Omega=s_{1}\Omega$}
\put(-70,90){$s_{2}\Omega$}
\put(-170,170){$s_{1}s_{2}\Omega$}
\put(-250,330){$s_{1}^{2}s_{2}\Omega$}
\put(-130,330){$s_{2}^{2}\Omega$}
\put(70,330){$s_{1}s_{3}\Omega$}
\put(190,330){$s_{4}\Omega$}
\put(30,210){$s_{3}\Omega$}
}
}

%%%%%%%%%%%%%%%%%%%%%%%%%%%%%%%%%%%%%%%%%
\def\bose{
\put(-250,-60){$P_{2}(1)|_{{\cal B}}=P_{\infty}(1)|_{{\cal B}}= $ Bose-Fock}
\put(0,0){\vertices}
\put(10,10){\line(0,1){100}}
\put(10,110){\line(1,1){100}}
\qbezier(10,10)(-90,10)(-93,210)
\qbezier(10,10)(-190,10)(-190,310)
\multiput(0,100)(-10,10){20}{$\cdot$}
\multiput(104,200)(0,10){10}{$\cdot$}
%\put(10,110){\line(-1,1){100}}
\put(-90,210){\line(2,1){200}}
\put(110,210){\line(1,1){100}}
%\put(110,210){\line(0,1){100}}
\put(-90,210){\line(0,1){100}}
%\put(-90,210){\line(-1,1){100}}
\put(10,110){\line(-1,2){100}}
{\small
\put(40,0){$\Omega$}
\put(30,90){$b_{1}^{*}\Omega$}
\put(-160,180){$b_{2}^{*}\Omega$}
\put(-250,330){$b_{3}^{*}\Omega$}
\put(-130,330){$b_{1}^{*}b_{2}^{*}\Omega$}
\put(70,330){$(b_{2}^{*})^{2}\Omega$}
\put(190,330){$(b_{1}^{*})^{3}\Omega$}
\put(120,170){$(b_{1}^{*})^{2}\Omega$}
}
}
%%%%%%%%%%%%%%%%%%%%%%%%%%%%%%%%%%%%%%%%%

\noindent
\thicklines
%\framebox{
\setlength{\unitlength}{.1mm}
\begin{picture}(1200,600)(-50,0)
\put(220,200){\otwo}
\put(850,200){\oni}
\end{picture}
%}

\noindent
\thicklines
%\framebox{
\setlength{\unitlength}{.1mm}
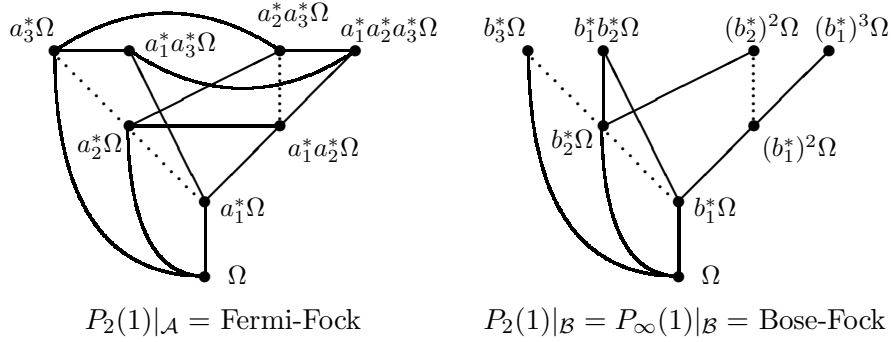
\begin{picture}(1200,500)(-50,0)
\put(220,100){\fermi}
\put(850,100){\bose}
\end{picture}
%}

\noindent
where $P_{2}(1)$ and $P_{\infty}(1)$ are as in $\S$ \ref{section:second},
vertices mean orthogonal vectors of the representation space of $\co{2}$
and these vectors are orthogonal each other
(see also graphs in $\S$ 3 of \cite{IWF01}).

%%%%%%%%%%%%%%%%%%%%%%%%%%%%%%%%%%%%%%%%%%%%%%%%%%
%
% subsection 4.2
%
\ssft{Computation of $U^{*}$ on the Fock representation}
\label{subsection:fourthtwo}
For $U$ in (\ref{eqn:unitary}), we show the adjoint operator $U^{*}$ of $U$.
Let ${\cal P}$ denote the set of all nonempty finite subsets of ${\bf N}$.
For $S\in {\cal P}$,
the {\it block decomposition} of $S$ is the partition 
$S=S_{1}\sqcup\cdots\sqcup S_{m}$
of $S$ such that there exists $x_{1},\ldots,x_{m}\in S$,
$k_{1},\ldots,k_{m}\in {\bf N}\cup\{0\}$ which satisfy
the following:
\begin{enumerate}
%(i)
\item
$S_{i}=\{x_{i},x_{i}+1,\ldots,x_{i}+k_{i}\}$ for each $i=1,\ldots,m$,
%(ii)
\item
$x_{i}+k_{i}+1<x_{i+1}$ for  each $i=1,\ldots,m-1$ when $m\geq 2$.
\end{enumerate}
For $S=\{n_{1},\ldots,n_{m}\}\in {\cal P}$, we write
$a_{S}^{*}\equiv a_{n_{1}}^{*}\cdots a_{n_{m}}^{*}$
when $n_{1}<\cdots <n_{m}$.
If $S=S_{1}\sqcup \cdots\sqcup S_{m}$ is the block decomposition of $S$, 
then $a_{S}^{*}=a_{S_{1}}^{*}\cdots a_{S_{m}}^{*}$.
From Theorem \ref{Thm:main}, we obtain the action of $U^{*}$. 
%
% Proposition 4.3
%
\begin{prop}
\label{prop:adjoint}
For $S\in {\cal P}$,
assume that $S=S_{1}\sqcup \cdots \sqcup S_{m}$
is the block decomposition of $S$ and 
$S_{i}\equiv\{n_{i},\ldots,n_{i}+l_{i}\}$ for $i=1,\ldots,m$.
Then the adjoint operator $U^{*}$ of $U$ in (\ref{eqn:unitary})
is given as follows:
\[
\begin{array}{rl}
U^{*}\Omega_{F}=&\Omega_{B},\\
\\
U^{*}a_{S}^{*}\Omega_{F}
=&
D\cdot (b_{n_{1}}^{*})^{l_{1}+1}
(b_{n_{2}-l_{1}-1}^{*})^{l_{2}+1}\cdots
(b_{n_{m}-\sum_{i=1}^{m-1}(l_{i}-1)}^{*})^{l_{m}+1}\Omega_{B}
\end{array}
\]
where $D$ denotes the normalization constant given by
\[D=\{(l_{1}+1)!\cdots (l_{m}+1)!\}^{-1/2}.\]
\end{prop}

%%%%%%%%%%%%%%%%%%%%%%%%%%%%%%%%%%%%%%%%%%%%%%%%%%%%%%
%
% subsection 4.3
%
\ssft{Realizations of Fock spaces as $\ltn$}
\label{subsection:fourththree}
We realize four irreducible representations $P_{2}(1)$,
$P_{\infty}(1)$ and Fock representations of bosons and fermions on
the Hilbert space $\ltn$.
Let $\{e_{n}:n\in {\bf N}\}$ denote the standard basis of $\ltn$.
Then
%
% Equation 4.1
%
\begin{equation}
\label{eqn:fourone}
\pi(t_{i})e_{n}\equiv e_{2(n-1)+i}\quad (i=1,2,\,n\in {\bf N})
\end{equation}
defines a representation of $\co{2}$ which is $P_{2}(1)$.
From Proposition \ref{prop:oni}(i), $\pi|_{\coni}$ is $P_{\infty}(1)$.
This is given as follows:
%
% Equation 4.2
%
\begin{equation}
\label{eqn:fourtwo}
\pi(s_{m})e_{n}=e_{2^{m-1}(2n-1)}\quad (m,n\in{\bf N}).
\end{equation}
We see that $\pi|_{{\cal A}}$ is the Fermi-Fock space with the vacuum $e_{1}$
(see also (3.40) of \cite{AK1}):
%
% Equation 4.3
%
\begin{equation}
\label{eqn:fourthree}
\pi(a_{m}^{*})e_{1}=e_{2^{m-1}+1},\quad\pi(a_{m})e_{1}=0\quad(m\geq 1).
\end{equation}
Furthermore,
we see that $\pi|_{{\cal B}}$ is the Bose-Fock space with the vacuum $e_{1}$:
%
% Equation 4.4
%
\begin{equation}
\label{eqn:fourfour}
\pi(b_{m}^{*})e_{1}=e_{2^{m-1}+1},\quad\pi(b_{m})e_{1}=0\quad(m\geq 1).
\end{equation}

%%%%%%%%%%%%%%%%%%%%%%%%%%%%%%%%%%%%%%%%%%%%%%%%%%%%

\end{document}